\documentclass[12pt]{article}

\def\c#1{{\cal #1}}
\def\1{{\bf 1}}
\newcommand{\eq}{\begin{equation}}
\newcommand{\en}{\end{equation}}
\newcommand{\eqa}{\begin{eqnarray}}
\newcommand{\ena}{\end{eqnarray}}

\begin{document}
\title{The hidden geometry of the quantum Euclidean space
\footnote{Talk given by the first author at the
``Quantum Group Symposium'', Prague, June 1998}}
\author{Gaetano Fiore$^{1,2}$, \  John Madore$^{3,4}$
        \and
        $\strut^1$Dip. di Matematica e Applicazioni, Fac.  di
Ingegneria\\ 
        Universit\`a di Napoli, V. Claudio 21, 80125 Napoli
        \and
        $\strut^2$I.N.F.N., Sezione di Napoli,\\
        Mostra d'Oltremare, Pad. 19, 80125 Napoli
        \and
        $\strut^3$Max-Planck-Institut f\"ur Physik 
        (Werner-Heisenberg-Institut)\\
        F\"ohringer Ring 6, D-80805 M\"unchen
        \and
        $\strut^4$Laboratoire de Physique Th\'eorique et Hautes
Energies\\
        Universit\'e de Paris-Sud, B\^atiment 211, F-91405 Orsay
        }

\setcounter{page}{1}   
\maketitle

\begin{abstract}
We briefly describe how to introduce the basic notions 
of noncommutative differential geometry on the 3-dim quantum
space covariant under the quantum group of rotations $SO_q(3)$.
\end{abstract}

\section{Introduction and preliminaries}

It is a rather old idea that the
micro-structure of space-time at the Planck level might be better
described using a noncommutative geometry. Here we consider the
formalism
of Dubois-Violette, Madore, Masson, Mourad, {\it et al.},
and apply it to the noncommutative algebra
describing the quantum Euclidean space ${\bf R}_q^3$
\cite{FRT}, namely the
quantum space covariant under the quantum group $SO_q(3)$.
This involves an interesting cross-fertilization
between the noncommutative-geometry formalism with the quantum space 
and
quantum group machinery. We briefly describe the main results of our
work
\cite{FioMad98'}. There,
we introduced a metric and an `almost' metric-compatible
linear connection on the quantum Euclidean 
space, equipped with its (two) standard $SO_q(3)$-covariant 
differential
calculi; correspondingly, the `frame' or dreibein has been also found. 
Modulo a conformal factor, which might however be
reabsorbed into a formulation of metric
compatibility more suitable for the present case, the curvature 
turns out to be zero, suggesting that the quantum space is flat as in
the 
commutative limit. In a separate
paper we shall show that in the same limit the traditional quantum 
space 
coordinates go to suitable general (non-cartesian) coordinates. 
This will allow a cure of some unpleasent features \cite{fiolat} of a
naive 
physical interpretation
of the representation theory of $Fun({\bf R}_q^3)$.

The preliminaries contained in this section are based  
especially on the works \cite{DubMadMasMou95, DimMad96};
for an introduction see Ref. \cite{Mad95}.
The starting point is a noncommutative algebra $\c{A}$ which has as
commutative limit 
the algebra of functions on some manifold $\c{M}$ and over
$\c{A}$ a differential calculus~\cite{Con94} 
$\{d,\Omega^*(\c{A})\}$ which
has as corresponding limit the ordinary de~Rham differential
calculus; as known~\cite{Con94}, $\{d,\Omega^*(\c{A})\}$ is completely
determined by the left and right module structure of the
$\c{A}$-module of 1-forms $\Omega^1(\c{A})$.  
By definition a {\it metric} is a $\c{A}$-bilinear map
\eq
g:\Omega^1(\c{A})\otimes_{\c{A}} \Omega^1(\c{A})\rightarrow
\c{A}.
\en
$\c{A}$-bilinearity means
\eq
g(f\xi\otimes\eta)=f g(\xi\otimes\eta),\qquad\qquad
g(\xi\otimes\eta f)=g(\xi\otimes\eta) f,
\en
for any $f\in\c{A}$ and $\xi, \eta\in\Omega^1(\c{A})$. This is
a definition in the ``cotangent space of the deformed manifold'';
one could also formulate it in the ``tangent space''.
In the commutative limit $\c{A}$-bilinearity is equivalent to the very
important
requirement of locality  of $g$ in both arguments at each point
$x\in\c{M}$:
\eq
[g(f\xi\otimes\eta)](x)=f(x)\cdot [g(\xi\otimes\eta)](x),
\qquad\qquad
[g(\xi\otimes\eta f)](x)=[g(\xi\otimes\eta)](x) f(x).
\en

A {\it linear connection} is a map (cfr. \cite{Kos60})
\eq
D:\Omega^1(\c{A}) \rightarrow\Omega^1(\c{A})\otimes_{\c{A}}
\Omega^1(\c{A})
\en
together with a ``generalized flip'' $\sigma$, i.e. a $\c{A}$-bilinear
map
\eq
\sigma:\Omega^1(\c{A})\otimes_{\c{A}} \Omega^1(\c{A})
\rightarrow\Omega^1(\c{A})
\otimes_{\c{A}} \Omega^1(\c{A})
\en
such that $D$ satisfies the left and right Leibniz rules
\eqa
D (f \xi)  & =&  df \otimes \xi + f D\xi                   
\label{2.2.2}\\
D(\xi f) &= &\sigma (\xi \otimes df) + (D\xi) f .         
\label{second}
\ena
Let $\pi$ be the projection
\eq
\pi:\Omega^1(\c{A})\otimes_{\c{A}} \Omega^1(\c{A})
\rightarrow\Omega^2(\c{A}).
\en
The {\it torsion} is the map $\Theta= d-\pi\circ D$. The connection $D$
is torsion-free iff
\eq
\pi\circ (\sigma+ \1)=0.
\en
One can naturally extend $D$ to higher tensor powers, e.g.
\eq
D_2(\xi\otimes\eta)=
D\xi\otimes\eta+\sigma_{12}(\xi\otimes D\eta),
\en
where we have introduced the tensor notation $\sigma_{12}=\sigma\otimes
\1$.
The metric-compatibility condition for $g,D$ reads $g_{23}\circ
D_2=d\circ g$.

The {\it curvature} 
$\mbox{Curv}:\Omega^1(\c{A})\rightarrow\Omega^2(\c{A})
\otimes_{\c{A}} \Omega^1(\c{A})$ is defined by
\eq
\mbox{Curv} = \pi_{12}\circ D_2\circ D.                \label{curv}
\en
It is always left $\c{A}$-linear, and right $\c{A}$-linear only in
certain
models; in general, right linearity is guaranteed only in the
commutative
limit. Therefore in this limit the curvature is local,  an 
essential physical requirement for a reasonable definition of a
curvature.

If $\c{A},\Omega^1(\c{A})$ are $*$-algebras and $d$ is real,
$(df)^*=df^*$,  $D$ can be made also real if we define \cite{FioMad98}
the involution
on $ \Omega^1(\c{A})\otimes_{\c{A}} \Omega^1(\c{A})$ by 
\eq
(\xi\otimes\eta)^*=\sigma(\eta\otimes\xi^*)
\label{star1}
\en
(note that this expression has the correct classical limit), with a
$\sigma$ such that the square of $*$ gives the identity. So
real structures on the
tensor product are in one-to-one correspondence with right Leibniz
rules.
$D_2$ is real iff $\sigma$ in addition fulfils the braid equation
\eq
\sigma_{12}\sigma_{23}\sigma_{12}=\sigma_{23}\sigma_{12}
\sigma_{23}.                                                 
\label{genbraid}
\en
The curvature is real if (\ref{star1}), (\ref{genbraid}) are satisfied.

Now assume that there exists a {\it frame}, i.e.
a special basis $\theta^a\in\Omega^1(\c{A})$,  $1\leq a \leq n$, such
that
\eq
[\theta^a,\c{A}]=0 
\en
and any $\xi\in\Omega^1(\c{A})$ can be uniquely written in the form
$\xi_a\theta^a$, with $\xi_a\in\c{A}$. This is possible only
if the limit manifold $\c{M}$ is parallelizable. It has the advantage 
that for any $f\in\c{A}$ the computation of
commutator $[\xi,f]$ is reduced to the computation of
the commutators $[\xi_a,f]$ in $\c{A}$.
Assume also that there exist $n$ inner derivations $e_a$, 
\eq
e_a f:=
[\lambda_a, f]
\en
($\lambda_a\in\c{A}$), dual to $\theta^a$: $\theta^a(e_b)=\delta^a_b$.
Then
\eq
\theta := - \lambda_a \theta^a                                
\label{dirac}
\en 
is the `Dirac operator'~\cite{Con94} for d:
\eq
df = -[\theta,f].                                             
\label{extra}
\en
$\theta^a$ is a very convenient basis to work with. For instance, from 
$\c{A}$-bilinearity it immediately follows that  the elements
\eq
g^{ab}:=g(\theta^a\otimes \theta^b)
\en
lie in the center $\c{Z}(\c{A})$ of $\c{A}$. We shall be interested in
the
case that $\c{Z}(\c{A})={\bf C}$. In the commutative limit
the condition $g^{ab}\in{\bf C}$ characterizes the vielbein or  
`moving frame' of E. Cartan, which is determined
up to a linear transformation; if this condition is fulfilled for any 
value of the deformation parameter the $\theta^a$ remain
uniquely determined up to a linear transformation
and are particularly convenient objects to be used to guess 
a physically sensible formulation of noncommutative-geometric notions. 

\section{Application of the formalism to the quantum Euclidean space}

Take  `the algebra of functions on the quantum Eucldean
space ${\bf R}_q^3$ 
\cite{FRT} as $\c{A}$ and over it one of the two $SO_q(3)$-covariant
differential calculi \cite{cawa}. The treatment
of the other calculus can be done in a completely
parallel way, see ref. \cite{FioMad98'}.
Here we are interested in the case of
a real positive $q$. We ask if they fit in the
previous scheme. We shall denote by $\Vert \hat R^{ij}_{hk}\Vert$ 
the braid matrix
of $SO_q(3)$, by $g_{ij}=g^{ij}$ the $SO_q(3)$-covariant
metric; here and below all indices will take the
values $-,0,+$. In the commutative limit $q\rightarrow 1$
$g_{ij}\rightarrow \delta^{i,-j}$.
The projector dedomposition of $\hat R$ is
\eq
\hat R = q\c{P}_s - q^{-1}\c{P}_a + q^{-2}\c{P}_t ;           
\label{decompo}
\en
$\c{P}_s$, $\c{P}_a$, $\c{P}_t$ are $SO_q(3)$-covariant
$q$-deformations of respectively the symmetric trace-free,
antisymmetric and trace projectors. The trace projector is
1-dimensional and is related to $g_{ij}$ by
\eq
\c{P}_t{}_{kl}^{ij} \propto g^{ij}g_{kl},         \label{bibi}
\en
$\c{A}$ is generated by $x^-,x^0,x^+$ fulfilling $\c{P}_a xx=0$,
or more explicitly
\eqa
&&x^- x^0 = q\, x^0 x^-,\nonumber\\
&&x^+ x^0 = q^{-1} x^0 x^+, \label{xxcr3}\\
&&[x^+, x^-] = h (x^0)^2.\nonumber
\ena
where we define $h = \sqrt q - 1/\sqrt q$. The real structure on $\c{A}$
is defined by $(x^i)^*=x^jg_{ji}$, or more explicitly
\eq
(x^-)^* = \sqrt q x^+, \qquad (x^0)^* = x^0, \qquad
(x^+)^* = 1/\sqrt{q} x^-.
\en
$\c{Z}(\c{A})$ is generated by the $SO_q(3)$-covariant real element
\eq
r^2 := g_{ij} x^i x^j = \sqrt{q}x^+x^-+ (x^0)^2+1/\sqrt{q}x^-x^+.
\label{sql}
\en
Let $\xi^i = dx^i$. One $SO_q(3)$-covariant calculus, which we shall
denoted by $\{d,\Omega^*(\c{A})\}$, is determined by the commutation
relations
\eq
x^i \xi^j = q\,\hat R^{ij}_{kl} \xi^k x^l.                    
\label{xxicr}
\en
Unfortunately neither calculus has a real exterior derivative, 
and up to now no way was known to make it closed under involution
\cite{olezu}; rather, each exterior algebra is mapped into
the other under the natural involution.
The `Dirac operator' (\ref{extra}) corresponding to $d$ is the
$SO_q(3)$-invariant element \cite{zumino}
$\theta := (q-1)^{-1} q^2 r^{-2} x^i\xi^j g_{ij}$;
note that $\theta$ is singular in the commutative limit.

In our work \cite{FioMad98'} we have found the following new results.
\begin{enumerate}

\item There exist two torsion-free, `almost' metric-compatible
linear connections, given by the formula
\eq
D_{(0)} \xi = -\theta \otimes \xi + \sigma_{0}
(\xi \otimes \theta)                            \label{conn}
\en 
The two corresponding generalized
flips $\sigma_{0}$ are determined by $S=q\hat R,(q\hat R)^{-1}$,
where $S$ is the $\c{A}$-valued matrix defined by
\eq
\sigma_0 (\xi^i\otimes \xi^j)  =: S^{ij}{}_{hk}\,\xi^h\otimes \xi^k;
\en
its knowledge allows one to extend by linearity $\sigma_{0}$ to all 
$\Omega^1(\c{A})\otimes_{\c{A}}\Omega^1(\c{A})$ in a unique way.
$D_{(0)}$ `almost' metric-compatible means compatible 
up to a conformal factor with the metric given in the next item;
a strict compatibility does not seem possible.
Both $\sigma_{0}$ fulfil the braid equation (\ref{genbraid})
and both $D_{(0)}$ are $SO_q(3)$-invariant.

\item If we extend $\c{A}$ by adding the `dilatation' generator 
$\Lambda$  
\eq
x^i\Lambda= q\Lambda x^i  
\en
together with its inverse $\Lambda^{-1}$ (we shall normalize them 
so that $\Lambda^*=\Lambda^{-1}$) and set $d\Lambda=0$, then 
up to normalization there exists a unique metric $g_0$, 
\eq
g_0(\xi^i\otimes \xi^j)  = g^{ij}\,r^2\Lambda^2
\en
($g_{ij}$ is the $SO_q(3)$-covariant metric matrix),
which is compatible with the two $D_{(0)}$ up to the conformal
factors $q^2,q^{-2}$,
\eq
S^{ij}{}_{hk}g^{kl}S^{mn}{}_{jl}=q^{\pm 2}g^{im}\delta^n_h,
\en
respectively in the cases $S=(q\hat R)^{\pm 1}$. A strict compatibility
would have required no $q^{\pm 2}$ at the rhs.

\item Curv=0 for both $D_{(0)}$.

\item If we further extend $\c{A}$ by adding also the generators
$r$ [the square root of (\ref{sql})], its inverse $r^{-1}$ and
the inverse $(x^0)^{-1}$ of $x^0$, then there exist a frame
$\theta^a$ , $a=-,0,+$,
and a dual basis $e_a$ of inner derivations given by 
\eq
\theta^a := \Lambda^{-1}\, \theta^a_i \xi^i            \label{ok}
\en
with
\eq
\Vert\theta^a_i\Vert:= \left\Vert
\begin{array}{ccc}
(x^0)^{-1}                            &                  & \cr
\sqrt{q}(q+1) (rx^0)^{-1}x^+          & r^{-1}           &\cr
-\sqrt{q}q(q+1) (r^2x^0)^{-1}(x^+)^2  & - (q+1)r^{-2}x^+ & r^{-2}x^0 
\cr
\end{array}
\right\Vert
\en
\eq
\begin{array}{l}
\lambda_- = + h^{-1} q \Lambda(x^0)^{-1} x^+,         \\   
\lambda_0 = - h^{-1} \sqrt q \Lambda (x^0)^{-1} r,\\
\lambda_+ = - h^{-1} \Lambda (x^0)^{-1} x^-.
\end{array}                                                    
\label{lambda'}
\en
$e_a x^i= q\Lambda \, e^i_a$, where $\Vert e^i_a\Vert$ is 
(left and right) inverse of the $\c{A}$-valued matrix
$\Vert\theta^a_i\Vert$.
Its elements fulfil the `$RTT$-relations' \cite{FRT}
\eq
\hat R_{kl}^{ij}\, e^k_a e^l_b
= e^i_c e^j_d \,\hat R_{ab}^{cd}                            \label{RTT}
\en
as well as the `$gTT$-relations'
\eq
g^{ab} e^i_a e^j_b = r^2 g^{ij}\qquad\qquad
g_{ij} e^i_a e^j_b = r^2 g_{ab}                             \label{gTT}.
\en
In a sense $r^{-1}e^i_a$ are a realization of the generators
$T^i_a$ of $SO_q(3)$.
As a consequence we find
\eq
\c{P}_t{}^{ab}_{cd} \theta^c \theta^d = 0 \qquad\qquad
\c{P}_s{}^{ab}_{cd} \theta^c \theta^d = 0,               \label{ththcr}
\en
the same commutation relations fulfilled by the $\xi^i$'s.
Finally, up to a normalization
$g_0(\theta^a\otimes\theta^b) = g^{ab}$ .

\item $\Omega^*(\c{A})$ is closed under the involution defined by
\eq
(x^i)^{\star}= x^jg_{ji} \qquad\qquad
(\theta^a)^{\star} =\theta^b g_{ba} 
\en
(the latter acts nonlinearly on the $\xi^i$'s: 
$(\xi^i)^{\star}=\Lambda^{-2} \xi^j\,c_{ji}$,
with $c_{ji}\in\c{A}$).

\end{enumerate}

The reality structure of these differential calculi is an old
but always present problem (see \cite{olezu}). The solution
proposed in item 5 is not fully satisfactory, at least naively.
For instance, it does not yield real $d,D$; only the curvature 
is real, for the simple reason
that it vanishes. The involution cannot be consistently extended
to $\Omega^*(\c{A})\otimes\Omega^*(\c{A})$ according
to (\ref{star1}). Finally, apparently it has not
the correct classical limit. 

A more careful analysis is needed at this point, but is out
of the scope of the present report (for more details
see Ref. \cite{FioMad98'}). It involves the investigation
of the properties of the $*$-representations of $\Omega^*(\c{A})$
and seems to suggest a
more sophisticated version of the proposal in item 5, in which
the opposite properties of the two differential calculi
cancel with each other. The problems mentioned above and the fact that
the
linear connections $D_{(0)}$ are metric-compatible up
to conformal factors (or, in other words, are only conformally flat)
may be related, in the sense that a satisfactory formulation
of the reality properties could eventually yield also a new and
satisfactory 
formulation of metric-compatibility which can be strictly fulfilled.
A careful analysis of the commutative limit is also needed 
in order to propose a reasonable correspondence principle between the
`new' theory and classical differential geometry.

\end{document}